\documentclass[12pt,reqno,notitlepage]{amsart}

\pdfoutput=1
\usepackage{amsmath,amsfonts,amsthm,amssymb}
\usepackage{enumerate}
\usepackage{indentfirst}
\usepackage[dvips]{graphicx}
\usepackage[all]{xy}
\usepackage{stackrel}
\usepackage{microtype}

\usepackage{marginnote}

\usepackage[latin1]{inputenc}
\usepackage[T1]{fontenc}
\usepackage{verbatim}
\usepackage{color}
\usepackage[normalem]{ulem}

\usepackage{tikz}

\usepackage{hyperref}  
\hypersetup{pdfborder={0 0 0}, 
colorlinks=true, 
citecolor=blue,
linktoc=page,
pdfauthor={}, 
pdftitle={}
}

\usepackage[centering, includeheadfoot, hmargin=1.0in, tmargin=1.0in, 
bmargin=1in, headheight=6pt]{geometry}
\usepackage{lmodern}

\newcommand{\Fcal}{\mathcal{F}}

\newcommand{\Z}{\mathbb{Z}}
\renewcommand{\P}{\mathbb{P}}

\newcommand{\C}{\mathbb{C}}

\DeclareMathOperator{\euler}{\chi_\mathrm{top}}

\newtheorem{theorem}{Theorem}[section]
\newtheorem{lemma}[theorem]{Lemma}
\newtheorem{corollary}[theorem]{Corollary}
\newtheorem{proposition}[theorem]{Proposition}

\theoremstyle{definition}

\newtheorem{remark}[theorem]{Remark}
\newtheorem{example}[theorem]{Example}

\theoremstyle{plain}

\newcommand{\terminou}{\hfill$\lrcorner$}

\begin{document}

\author[{\scriptsize\rm{}T.~Fassarella and N.~Medeiros}]
{
T. Fassarella \and N. Medeiros 
}
\title[{\scriptsize\rm{}Arrangements and the Gauss map of a pencil}]
{Hyperplane arrangements and the Gauss map of a pencil}

\subjclass[2020]{Primary 14E05, 05B35; Secondary 14B05, 37F75.}
\keywords{{Hyperplane} arrangements,  characteristic polynomial, Gauss map, multidegrees}
\thanks{
The authors are supported by CNPq Universal 10/2023. 
The first author is supported by CAPES-COFECUB programme project number: Ma 1017/24.
}

\begin{abstract} 
We show that the coefficients of the characteristic polynomial of a central 
{affine complex} 
hyperplane arrangement $\mathcal A$, coincide with the multidegrees of the Gauss map of a pencil of hypersurfaces naturally associated to 
$\mathcal A$. 
{As a consequence, we obtain a proof of the Heron-Rota-Welsh conjecture for matroids representable over a field of characteristic zero.} 
\end{abstract}
\maketitle



\section{Introduction}
This short paper {is devoted to the study of} 
the log-concavity and unimodality of the sequence of  coefficients of the characteristic polynomial $\chi_\mathcal A(t)$ associated with a central affine arrangement $\mathcal A$ over the field of complex numbers $\C$. A sequence of real numbers $a_0, \dots, a_{n}$ is  {\it log-concave}, if  for all $i\in\{1,\dots, n-1\}$, the following inequality holds
\[
a_i^2 \geq a_{i-1}a_{i+1}.
\]
It is called {\it unimodal} if there exists an index $i\in\{0,\dots,n\}$ such that 
\[
a_0\le \cdots \le a_{i-1}\le a_i \ge a_{i+1}\ge\cdots \ge a_{n}. 
\]

 The Heron-Rota-Welsh conjecture \cite{Rot70, Her72, Wel76}, for matroids representable over a field of characteristic zero, can be reduced to showing that the  sequence formed by the coefficients of $\chi_\mathcal A(t)$  is log-concave. This conjecture generalizes those of Read \cite{Re68} and Hoggar \cite{Hog74}, extending their scope to matroids: Read's conjecture posits that the absolute values of the coefficients of a chromatic polynomial of a graph form a unimodal sequence, while Hoggar's 
conjecture asserts that this sequence is log-concave.


J. Huh \cite[Corollary 27]{Huh12} proved the Read-Hoggar conjecture and the Heron-Rota-Welsh conjecture for matroids representable over a field of characteristic zero.  
{The proof is quite sophisticated}: it involves works of singularity theory, algebraic geometry, and convex geometry; it relies on a theorem of Dimca and Papadima \cite{DP03} and three technical lemmas from commutative algebra.
In his proof, J. Huh connected the coefficients of the characteristic polynomial with the multidegrees of the polar map associated with a hyperplane arrangement that naturally emerges from the matroid.  The present work is motivated by this connection: we show, in our Theorem~\ref{thm:main} below, that the coefficients of the characteristic polynomial are precisely the multidegrees of the Gauss map of a {certain} pencil of hypersurfaces associated to the arrangement. The log-concavity follows immediately from this. 
Our approach is based on the results of \cite{FP07}, which hinges on logarithmic foliations. 


Before going into specifics, we remark that K. Adiprasito, J. Huh, and E. Katz \cite{AHK18} proved the Heron-Rota-Welsh conjecture for matroids in full generality. For an excellent overview of the key developments leading up to their final result, we refer to the introduction of \cite{AHK18}.

Let $\mathcal A$ be a central hyperplane arrangement in $\C^{n+1}$, consisting of a collection of $k\ge 1$ subspaces of dimension $n$. 
Let $\chi_\mathcal A(t)$ be its characteristic polynomial, as defined in Section~\ref{sec:chararrang}.
The central arrangement $\mathcal A$ defines a projective arrangement, whose union
of hyperplanes may be written as $V(f)\subset \P^{n}$,  where $f$ is a product of $k$ homogeneous linear forms in the variables $x_0, \dots, x_n$.

By introducing an additional variable $x_{n+1}$, we can associate to $\mathcal A$ a pencil $\mathcal F$ of hypersurfaces in $\P^{n+1}$,  generated by $f$ and $x_{n+1}^k$. The members of {this pencil} are the hypersurfaces {$V_{a,b}$} defined by the  equation $af+bx_{n+1}^k$, where $(a:b)\in \P^1$. Now, consider the Gauss map $\mathcal G(\mathcal F): \P^{n+1} \dashrightarrow \check{\P}^{n+1}$
which sends a smooth point $p$  of a member $V_{a,b}$
to its tangent space $T_pV_{a,b}$ in the dual space $\check{\P}^{n+1}$.
This leads us to a sequence of integers, called {\it multidegrees}, associated to this rational map, see Section~\ref{sec:logfol} {or \cite[Section 7.1.3]{Dol24} for further details.}

\begin{theorem}
\label{thm:main}
With notations as above, let $d_0, \dotsc, d_{n+1}$ be the sequence of multidegrees of $\mathcal G(\mathcal F)$. Then
\[
\chi_\mathcal A(t) = d_0 t^{n+1}-d_1 t^n + \cdots + (-1)^{n+1} d_{n+1}. 
\]
\end{theorem}
\medskip

{In the main text, this is Theorem \ref{thm:main2}.}
\smallskip

It is straightforward to check that the sequence $d_0, \dots, d_{n+1}$ of multidegrees has no internal zeros, see Proposition~\ref{prop:nointernal}. Moreover, it follows from the Hodge-Khovanskii-Teissier type inequalities \cite[Example {1.6.4}]{Laz04} that this sequence is log-concave. 
See also \cite[Theorem 21]{Huh12} for a different proof and  a  complete characterization of representable homology classes of a product of two projective spaces.   Then the next result follows from  Theorem~\ref{thm:main}.

\begin{corollary}\label{cor:huhmatroid}{\rm (\cite{Huh12})}
Let $M$ be a matroid {representable over $\C$}. Then  the coefficients of the characteristic polynomial $\chi_M(t)$ form a log-concave sequence of integers with no internal zeros.
\end{corollary}

As mentioned earlier, Corollary~\ref{cor:huhmatroid} was previously conjectured  by Heron, Rota and Welsh. {It is worth noting that a matroid which is representable over a field of characteristic zero is also representable over $\C$, see the proof of \cite[Corollary 27]{Huh12} for further details.}  In the same context,  the chromatic polynomial $\chi_G(t)$ of 
{a finite} 
graph $G$, coincides with the characteristic polynomial of the graph arrangement $\mathcal A$ associated to $G$. Hence, since any  sequence of  nonnegative numbers that is log-concave and has no internal zeros is unimodal,  Theorem~\ref{thm:main} yields the following result, which was previously conjectured by Read {and Hoggar}.

\begin{corollary}\label{cor:huhgraph}{\rm (\cite{Huh12})}
{The coefficients of the chromatic polynomial $\chi_G(t)$ of any 
{finite}	
graph $G$ form a  {log-concave} sequence
{with no internal zeros}, 
hence their absolute values form a unimodal sequence.
}
\end{corollary}

\section{Characteristic polynomial of an arrangement}
\label{sec:chararrang}

We refer to \cite{Sta07} and \cite{OT92} for general background on hyperplane arrangements. 

Let $\mathcal A$  be 
a central
arrangement given by  a finite collection of affine hyperplanes  in $\C^{n+1}$.
We denote by $L_\mathcal A$ the set of all nonempty intersections of hyperplanes of $\mathcal A$,  with the ambient space $\C^{n+1}$ included as one of its elements.   Let us endow $L_\mathcal A$ with a structure of \emph{poset} (partially ordered set): define  $C_1 \le C_2$ in  $L_\mathcal A$ if $C_2\subseteq C_1$. Note that $\C^{n+1}\le C$ for all $C\in L_\mathcal A$. 

Let ${\rm Int}(L_\mathcal A)$ be the set of all closed intervals $[C_1, C_2]$ of $L_\mathcal A$. In order to define the characteristic polynomial of the arrangement, we need to introduce an important tool, the {\it M\"obius function}. It is a function  $\mu \colon {\rm Int}( L_\mathcal A) \to \Z$, defined by the following conditions:
\begin{itemize}
\item $\mu(C,C) = 1$ for all $C\in  L_\mathcal A $. 
\item $\mu(C_1,C_2) = -\sum_{C_1\le C < C_2} \mu(C_1,C)$ for all $C_1 < C_2$ in $ L_\mathcal A$
\end{itemize}
We write $\mu(C)=\mu(\C^{n+1}, C)$ for all $C\in L_\mathcal A$. 
The \emph{characteristic polynomial} $\chi_\mathcal A (t)$ of the arrangement is defined by
\[
\chi_\mathcal A (t) = \sum_{C\in L_{\mathcal A}} \mu(C)\cdot t^{\dim C}. 
\]

Given a hyperplane $H$ in $\mathcal{A}$, we define a triple of arrangements
$(\mathcal A, \mathcal A^*, \mathcal A^{**})$, where  $\mathcal A^* = \mathcal A - H$ is obtained from $\mathcal A$ by deleting $H$ and $\mathcal A^{**}$ is the arrangement in $H\simeq \C^n$ given by the restriction $\mathcal A^*|_H$. The characteristic polynomials of such a triple satisfy the fundamental deletion-restriction principle:
\begin{equation}\label{eq:chiA}
\chi_\mathcal A(t) = \chi_{\mathcal A^*}(t) - \chi_{\mathcal A^{**}}(t).
\end{equation}
See \cite[Lemma 2.2]{Sta07}.

Closely related to the characteristic polynomial, the \emph{chromatic polynomial} $\chi_G(t)$ of a graph $G$ counts the number of proper colorings using $t$ colors. A \emph{proper coloring} assigns colors to the vertices such that adjacent ones receive different colors. For a simple graph $G$ with vertex set $\{v_0, \dots, v_n\}$, there exists an associated arrangement $\mathcal A_G$ in $\C^{n+1}$ that satisfies
\begin{equation}
\label{eq:chiAchiG}
\chi_G(t)=\chi_{\mathcal A_G}(t).
 \end{equation}
This arrangement is defined by the union of hyperplanes $V(x_i-x_j)$, $i<j$, for all adjacent vertices $v_i$ and  $v_j$.  It turns out that $\chi_G(t)$ satisfies, similarly, the deletion-contraction principle for graphs 
\begin{equation*}
\label{eq:chiGprinciple}
\chi_G(t) = \chi_{G^*}(t) - \chi_{G^{**}}(t).
\end{equation*}
One approach to proving identity \eqref{eq:chiAchiG} is by induction on the number of edges of $G$, then on the number of hyperplanes of $\mathcal A_G$. See \cite[Theorem 2.7]{Sta07}. In Section \ref{sec:pencilarrang}, we will apply the same method to prove that $\chi_\mathcal A(t)$ coincides with the polynomial $\chi_\mathcal F(t)$, where the absolute value of each  coefficient is a multidegree of the Gauss map of a pencil associated to the arrangement.

\section{Logarithmic foliations on projective spaces} 
\label{sec:logfol}

A {\it logarithmic foliation} $\mathcal F$ on $\P^n$ is induced by  a rational 1-form
\begin{eqnarray}
\label{omegalog}
\omega = \sum_{i=1}^s \lambda_i\frac{df_i}{f_i}
\end{eqnarray}
where $f_1, \dots, f_s$ are irreducible homogeneous polynomials in $\C[x_0, \dots, x_n]$ and  $\lambda_1, \dotsc, \lambda_s$ are nonzero complex numbers satisfying the relation
\[
\sum_{i=1}^s \lambda_i \deg (f_i) = 0. 
\]
This condition ensures that the contraction $i_R\omega$ with the radial vector field $R=\sum_{i=0}^nx_i \frac{\partial}{\partial x_i}$  vanishes identically. Furthermore,  $\omega$ is {\it integrable}, in the sense that  
\[
d\omega \wedge \omega = 0.
\]

The $1$-form $\omega$ defines a global section of the sheaf $\Omega^1_{\P^n}(\log D)$, whose sections are 1-forms with at most simple poles along the divisor 
\[
D = V(f_1\cdots f_s) \subset\P^n
\]
called the \emph{polar locus} of $\omega$.
In a sufficiently small analytic neighborhood of a nonsingular point of $\omega$,   the integrability condition gives rise to a  holomorphic fibration whose relative tangent sheaf coincides with the subsheaf  of $T\P^n$ determined by the kernel of $\omega$. The analytic continuation of a fiber of this local fibration defines a  {\it leaf} of $\mathcal F$. 

At a nonsingular point $p\in \P^n$
of $\omega$, we denote by $T_p\mathcal F$ the {\it tangent space} of the leaf of $\mathcal F$ passing through $p$.  The assignment $p \mapsto T_p\mathcal F$ determines a rational map
\[
\mathcal G(\Fcal)\colon \P^n \dashrightarrow \check{\P}^n
\] 
called the \emph{Gauss map} of the foliation $\mathcal F$. Here, $\check{\P}^n$ denotes the dual space, formed by hyperplanes in $\P^n$.  Using the natural identification $\check{\P}^n\simeq \P^n$, we can see the Gauss map as a rational map $\mathcal G(\Fcal)\colon \P^n \dashrightarrow \P^n$. Moreover,  if we write $\omega$ as 
\[
\omega = \sum_{i=0}^{n} a_i dx_i
\] 
for suitable rational functions $a_0, \dots, a_n$, then $\mathcal G(\Fcal)$ is given by  $( a_0: \dots : a_n )$. 

We associate to $\mathcal F$ a sequence of integers, called {\it multidegrees}, as follows. If $\mathbb P^{n-j}\subset \P^{n}$ is a general linear  subspace of codimension $0\leq j \leq n$, then set
$d_j(\mathcal G(\Fcal))$, or simply $d_j$, {to be} the degree of the closed subset $\overline{\mathcal G(\Fcal)^{-1}(\P^{n-j})}$. 
In particular, 
\[
 d_0=1\quad  \text{and} \quad d_{n} = \deg (\mathcal G(\Fcal))
\]
where $\deg (\mathcal G(\Fcal))$ is the topological degree of the rational map $\mathcal G(\Fcal)$.

In order to apply some results of \cite{FP07}, we need to introduce some notation. Let $\omega$ be a rational logarithmic 1-form as in \eqref{omegalog}. 
Let $\pi\colon X \to \P^n$ be a sequence of blowups, along smooth centers of the polar locus $D$, such that $\pi^*(D)$ has only normal crossings.  We can see $\pi^*(\omega)$ as a global section of $\Omega^1_X(\log( \pi^*(D)))$.  Given an irreducible component $E$ of $\pi^*(D)$, 
one can {locally express} $\pi^*(\omega)$ at a general point of $E$ as 
\[
\pi^*(\omega) = \lambda(E)\cdot \frac{dh}{h} + \eta
\]
where $\lambda(E)\in \mathbb C$, $h$ is a local equation for $E$ and $\eta$ is a holomorphic 1-form.  We say that $\lambda(E)$ is the {\it residue} of $\pi^*(\omega)$ with respect to $E$. It is known that there are nonnegative integers $m_1, \dots, m_s$ such that 
\[
\lambda(E) = \sum_{i=1}^{s} m_i\lambda_i. 
\]
See for example \cite[Lemma 3]{FP07}. We say that $(\lambda_1, \dots, \lambda_s)$ is {\it nonresonant}, with respect to $\pi$, if $\lambda(E)\neq 0$ for every irreducible component $E$ of $\pi^*(D)$. 

Let $\mathcal F|_{\P^i}$ be the restriction of the foliation $\mathcal F$ to a general linear subspace $\P^i\subset \P^n$ of dimension $i$. We note that $\mathcal F|_{\P^i}$ is still a logarithmic foliation, defined by the 1-form $\omega|_{\P^i}$. In the next result, we describe the topological degree of the Gauss map $\mathcal G(\mathcal F|_{\P^i})$ of this foliation.  

\begin{lemma}
\label{lemma:versionFP}
Let the notation be as above. Assume that $(\lambda_1, \dots, \lambda_s)$ is nonresonant. Let $\P^i\subset\P^n$ denote a general linear space of dimension $i$ and let $\mathcal G(\mathcal F|_{\P^i})$ be the Gauss map of the restriction $\mathcal F|_{\P^i}$. Then, for $i=2, \dots, n$, we have
\[
d_i(\mathcal G(\mathcal F|_{\P^i})) = (-1)^{i-1} \euler (\P^{i-1}\backslash D). 
\]
\end{lemma}

\proof
Since $d_i(\mathcal G(\mathcal F|_{\P^i}))$ is the topological degree of $\mathcal G(\mathcal F|_{\P^i})$, the result follows from \cite[Proposition 1]{FP07} and the logarithmic Gauss-Bonnet theorem  ({see} \cite{Nor78}, \cite{Sil96} or \cite[p.\thinspace{}262]{Kaw78}). 
\endproof

We {will} apply the above lemma as follows. Let $f$ be a nonconstant homogeneous polynomial in $\C[x_0, \dots, x_n]$, of degree $k$. By introducing an additional variable $x_{n+1}$, 
we associate to $f$ a logarithmic foliation $\mathcal F$ on $\P^{n+1}$ given by
\begin{equation}
\label{eq:1form}
\omega = \frac{df}{f} - k \frac{dx_{n+1}}{x_{n+1}}. 
\end{equation}
We can see $\mathcal F$ as a pencil of hypersurfaces in $\P^{n+1}$, generated by $f$ and $x_{n+1}^k$. {Its} Gauss map $\mathcal G(\Fcal)$
sends a general point $p\in \P^{n+1}$ to the tangent space of the member of $\mathcal F$ passing through $p$. 

\begin{proposition}
\label{prop:versionFP}
Let $\mathcal F {\subset\P^{n+1}}$ be the logarithmic foliation associated to a nonconstant polynomial
$f\in \C[x_0, \dots, x_n]$ of degree $k$,
defined by the rational 1-form \eqref{eq:1form}.
Let $\mathcal G(\Fcal)$ be its Gauss map. Then, for $i=2, \dotsc, n+1$, we have
\[
d_i(\mathcal G(\mathcal F|_{\P^i})) = (-1)^{i-1} \euler (\P^{i-1}\backslash (V(f)\cup \P^{i-2}))
\]
where $\P^{i-2}\subset \P^{i-1}\subset \P^i$ are general linear spaces.
\end{proposition}

\proof
Let $V_i = V(f)\cap \P^i$ and define {$\P^{i-1} = V(x_{n+1})\cap\P^i$}, for any $i$.  Consider a resolution  $\pi\colon  X \to \P^i$ of $V_i$, that is, a sequence of blowups  such that the total transform $\pi^*(V_i)$ has only normal crossings. This process also yields a resolution of  $V_i\cup \P^{i-1}$. Now, 
write $f = f_1^{n_1}\cdots f_s^{n_s}$ as product of irreducible factors. Then $\mathcal F$ is defined by
\[
\omega = \sum_{i=1}^{s}n_i\frac{df_i}{f_i} - k \frac{dx_{n+1}}{x_{n+1}}. 
\]
We claim that the vector of residues $(n_1, \dots, n_s, -k)$ is nonresonant. Indeed, $\pi^*(\P^{i-1})$ has a unique irreducible component, with residue $-k$, because no center of the  blowups lies in $\P^{i-1}$.  Besides this, any other irreducible component of $\pi^*(V_i)$ has residue of the form $\sum_{i=1}^{s} m_in_i$,  for suitable nonnegative integers $m_1, \dotsc, m_s$; see \cite[Lemma 3]{FP07}. 

The conclusion now follows from Lemma~\ref{lemma:versionFP} by taking $D = V(f\cdot x_{n+1})$ and $\P^{i-2} = V(x_{n+1})\cap \P^{i-1}$.  
\endproof

\section{The pencil of an arrangement}
\label{sec:pencilarrang}

We keep the notations of Sections~\ref{sec:chararrang} and \ref{sec:logfol}. We associate to a central arrangement $\mathcal A$ in $\C^{n+1}$ the pencil $\mathcal F$ of hypersurfaces in $\P^{n+1}$ generated by $f(x_0,\dotsc,x_n)$ and $x_{n+1}^k$, where $f = h_1\cdots h_k$, the product of the homogeneous linear forms defining the members of the arrangement. 
Then $\mathcal F$ has $V(af+bx_{n+1}^k)$, for $(a:b)\in \P^1$, as its members. The pencil $\mathcal F$ can be seen as a logarithmic foliation defined by the rational 1-form
\[
\frac{{d} f}{f} - k \frac{{d} x_{n+1}}{x_{n+1}}. 
\]
Let $\mathcal G(\mathcal F):\P^{n+1} \dashrightarrow \P^{n+1}$ be the Gauss map of $\mathcal F$, given by the linear system 
\[
\langle x_{n+1}\frac{\partial f}{\partial x_0}, \dotsc, x_{n+1}\frac{\partial f}{\partial x_n}, -kf \rangle.
\] 
As in Section~\ref{sec:logfol}, we consider a sequence of integers $d_0(\mathcal G(\mathcal F)), \dotsc, d_{n+1}(\mathcal G(\mathcal F))$, or simply, $d_0, \dots, d_{n+1}$, called  multidegrees of the rational map $\mathcal G(\mathcal F)$.
They satisfy 
\[
 d_0=1,\quad {d_1=k} \quad \text{and} \quad d_{n+1} = \deg (\mathcal G(\mathcal F))
\]
where $\deg (\mathcal G(\mathcal F))$ is the topological degree of  $\mathcal G(\mathcal F)$. The following result is an easy and useful fact, which holds for any rational map.

\begin{proposition}\label{prop:nointernal}
The above sequence  $d_0, \dots, d_{n+1}$ has no internal zeros. 
\end{proposition}
\proof
Since $d_0=1$, we assume that $d_j> 0$ for some $j\ge 1$. This implies that the dimension of the image $\operatorname{Im} \mathcal G(\mathcal F)$ of the rational map $\mathcal G(\mathcal F)$ is at least $j$. Then the intersection between $\operatorname{Im} \mathcal G(\mathcal F)$ and a general $\P^{n+1-i}$ is nonempty for all $i\le j$. This gives $d_i>0$ for all $i\le j$. 
\endproof

We consider the following polynomial defined from the sequence of multidegrees 
\[
\chi_\mathcal F (t) = d_0 t^{n+1}-d_1t^n+ \cdots + (-1)^{n+1} d_{n+1}.
\]
Note that for an arrangement with a single hyperplane we have 
\begin{equation}\label{eq:inicial}
\chi_\mathcal A (t) =  \chi_\mathcal F (t) = t^{n+1} - t^n. 
\end{equation}

Let us setup a bit more of notation.
Given $H\in \mathcal A$, denote by $\P^{n-1}_H\subset \P^n$ the corresponding hyperplane in $\P^n$.
As we have done in Section~\ref{sec:chararrang} for affine arrangements, we consider a
triple $(A,A^{*},A^{**})$ of projective arrangements, where $A, A^* \subset \P^n$ are the union of hyperplanes in $\mathcal A, \mathcal A^*$ respectively, and $A^{**} = A^*\cap \P^{n-1}_H$. 
Let $\mathcal F^*$ and $\mathcal F^{**}$ be the corresponding pencils associated to $\mathcal A^*$ and $\mathcal A^{**}$, respectively.

{Now we arrive at the central result of this note: The deletion-restriction principle holds for the polynomial $\chi_\mathcal F(t)$.}
\begin{proposition}
\label{prop:deletion}
With above notations, we have
\begin{equation}\label{eq:chiF}
\chi_\mathcal F(t) = \chi_{\mathcal F^*}(t) - \chi_{\mathcal F^{**}}(t).
\end{equation}
\end{proposition}
\proof
This is equivalent to show that the corresponding multidegrees satisfy the identity:
\begin{equation}\label{eq:multidegF}
d_i(\mathcal G(\mathcal F)) = d_i(\mathcal G(\mathcal F^*)) + d_{i-1}(\mathcal G(\mathcal F^{**}))
\end{equation}
for $i=1, \dotsc, n+1$. Since $d_1$ is the number of hyperplanes of the arrangement and $d_0=1$,  the identity holds for $i=1$.

Let us start with the case $i=n+1$. Our strategy is to use the inclusion-exclusion principle for the {topological} Euler characteristic.
We have 
\begin{equation}
	\label{uniondisjoint}
	\begin{aligned}
	A\cup \P^{n-1} =&\  A^*\cup \P^{n-1}_H \cup \P^{n-1} \\
	=&\  (A^*\cup \P^{n-1}) \cup [\P^{n-1}_H\setminus (A^* \cup \P^{n-1})]\\
	=&\  (A^*\cup \P^{n-1}) \cup [\P^{n-1}_H\setminus (A^{**}\cup \P^{n-2})]
\end{aligned}
\end{equation}
where $\P^{n-1}$ is a general hyperplane and $\P^{n-2}=\P^{n-1}\cap \P^{n-1}_H$. Since this is a disjoint union, we have
\[
(-1)^n\euler (\P^n\setminus (A\cup \P^{n-1})) = (-1)^n \euler (\P^n\setminus(A^*\cup \P^{n-1})) + (-1)^{n-1} \euler(\P^{n-1}_H\setminus (A^{**}\cup \P^{n-2}))
\]
hence by Proposition \ref{prop:versionFP} we obtain
\begin{equation*}	
	\label{eq-dn+1}
	d_{n+1}(\mathcal G(\mathcal F)) = d_{n+1}(\mathcal G(\mathcal F^*))  + d_n(\mathcal G(\mathcal F^{**}))
\end{equation*}
therefore concluding that \eqref{eq:multidegF} holds for $i=n+1$. 

Now let us assume that
$i\in\{2,\dots, n\}$. The argument in this case goes in a similar fashion. The key result here is \cite[Theorem 1]{FP07}, which
express the degree $d_i(\mathcal G(\mathcal F))$ in terms of degrees of restrictions to $\P^i$ and $\P^{i+1}$, namely
\begin{equation}
\label{eq:diF}
d_i(\mathcal G(\mathcal F)) = d_i(\mathcal G(\mathcal F|_{\P^i})) + d_{i+1}(\mathcal G(\mathcal F|_{\P^{i+1}})).
\end{equation}
Our task is to evaluate the terms on the right hand side. By Proposition~\ref{prop:versionFP} we have
\begin{equation*}\label{eq:eulerjvp2}
	d_{i}(\mathcal G(\mathcal F|_{\P^i})) = (-1)^{i-1} \euler (\P^{i-1}\setminus (A_{i-1}\cup \P^{i-2}))
\end{equation*}
where $A_{i-1} = A\cap \P^{i-1}$. Now, let $\P^i_H = \P^{n-1}_H\cap \P^i$,  $A_i^*=A^*\cap \P^i$ and $A_i^{**} = A^{**}\cap \P^i$.  By restricting to $\P^{i-1}$, in the same manner as we did in \eqref{uniondisjoint}, we can write $A_{i-1}\cup \P^{i-2}$ as a disjoint union
	\begin{eqnarray*}\label{eq:restrictionagain}
		A_{i-1}\cup \P^{i-2} =  (A_{i-1}^*\cup \P^{i-2}) \cup [\P^{i-2}_H\setminus (A_{i-2}^{**}\cup \P^{i-3})] \label{uniondisjoint2}
	\end{eqnarray*}
and from that we get
	\begin{eqnarray*}
		(-1)^{i-1}\euler (\P^{i-1}\setminus (A_{i-1}\cup \P^{i-2})) &=& (-1)^{i-1} \euler (\P^{i-1}\setminus(A_{i-1}^*\cup \P^{i-2})) \\
		&+& (-1)^{i-2} \euler(\P^{i-2}_H\setminus (A_{i-2}^{**}\cup \P^{i-3}))
	\end{eqnarray*}
where we agree that $A_{i-2}^{**}$ and $\P^{i-3}$ are empty if $i=2$. 
Again by Proposition \ref{prop:versionFP},
\begin{equation*}\label{eq:diFs}
	d_i(\mathcal G(\mathcal F|_{\P^i})) = d_i(\mathcal G(\mathcal F^*|_{\P^i})) + d_{i-1}(\mathcal G(\mathcal F^{**}|_{\P^{i-1}})).
\end{equation*}
In the same way we obtain the analogous expression for $d_{i+1}(\mathcal G(\mathcal F|_{\P^{i+1}}))$. 
By substituting into \eqref{eq:diF}, we have
\[
\begin{aligned}
d_i(\mathcal G(\mathcal F)) = & \  [d_i(\mathcal G(\mathcal F^*|_{\P^i})) + d_{i-1}(\mathcal G(\mathcal F^{**}|_{\P^{i-1}})) ] + [d_{i+1}(\mathcal G(\mathcal F^{*}|_{\P^{i+1}})) + d_{i}(\mathcal G(\mathcal F^{**}|_{\P^{i}}))] \\
= & \ [d_i(\mathcal G(\mathcal F^*|_{\P^i})) + d_{i+1}(\mathcal G(\mathcal F^{*}|_{\P^{i+1}})) ] + 
[  d_{i-1}(\mathcal G(\mathcal F^{**}|_{\P^{i-1}})) + d_{i}(\mathcal G(\mathcal F^{**}|_{\P^{i}}))] \\
= & \ d_i(\mathcal G(\mathcal F^*)) + d_{i-1}(\mathcal G(\mathcal F^{**}))
\end{aligned}	
\]
where, in the last equality, we applied \eqref{eq:diF} for $\mathcal F^*$ and $\mathcal F^{**}$. 
Therefore \eqref{eq:multidegF} holds, and that concludes the proof.
\endproof

The following result is Theorem~\ref{thm:main} of the introduction. 

\begin{theorem}
\label{thm:main2}
Let $d_0, \dots, d_{n+1}$ be the sequence of multidegrees of the Gauss map $\mathcal G(\mathcal F)$ associated to the arrangement $\mathcal A$ in $\C^{n+1}$.   Then 
\[
\chi_\mathcal A(t) = d_0 t^{n+1}-d_1t^n + \cdots + (-1)^{n+1} d_{n+1}. 
\]
\end{theorem}

\proof
The same deletion-restriction principle of Proposition~\ref{prop:deletion}  holds for $\chi_\mathcal A(t)$, see \eqref{eq:chiA}. Since the initialization process is assured by \eqref{eq:inicial}, the proof of the theorem follows by induction on the number of hyperplanes of the arrangement $\mathcal A$.  
\endproof

We emphasize the geometric meaning of the coefficients of $\chi_\mathcal A(t)$ when interpreted as multidegrees of the Gauss map.

\begin{remark}\label{rmk:critical}
Let $U$ denote the complement of $V(f\cdot x_{n+1})$ in $\P^{n+1}$ and let $\varphi$ be the rational function on $\P^{n+1}$ defined by $f/x_{n+1}^k$. We fix $a_0=1$ and for each $i=1,\dotsc, n$, let $a_i$ represent  the number of {\it critical points} of the rational function $\varphi|_{\P^i}$ in $U$. These critical points correspond to the  singularities  of the logarithmic 1-form $d\log (\varphi|_{\P^i})$, where $\P^i\subset \P^{n+1}$ is a general subspace of dimension $i$. The topological degree $d_{n+1}$ of the Gauss map $\mathcal G({\mathcal F})$ equals the total number of tangencies between a general $\P^n\subset \P^{n+1}$ and the members of $\mathcal F$. This gives 
\begin{equation}
d_{n+1} = a_n.
\end{equation}
Similarly, for the topological degree $d_i(\mathcal G(\mathcal F|_{\P^i}))$ of the Gauss map  of the restriction $\mathcal F|_{\P^i}$ to a general $\P^i$, we have $d_i(\mathcal G(\mathcal F|_{\P^i})) = a_{i-1}$. In particular, the identity \eqref{eq:diF} implies
\begin{equation}
d_i = a_{i-1}+a_i
\end{equation}
for all $i=1,\dotsc, n$. Consequently,
\begin{equation}\label{eq:critical}
\chi_\mathcal A(t) = (t-1)(a_0t^n-a_1t^{n-1}+\cdots+ (-1)^na_n). 
\end{equation}
It turns out that $a_i$ is the $i$-th multidegree of the gradient map ${\rm grad}f\colon\P^n \dashrightarrow \P^n$, see \cite[Corollary 2]{FP07}. Then \eqref{eq:critical} recovers one of the two identities of \cite[Corollary 25]{Huh12}.
\terminou
\end{remark}

\begin{example}
{Here is a basic example, taken from }\cite[Lecture 1]{Sta07}.  Figure \ref{fig:posetmobius2} shows the values of the M\"obius function for the arrangement in $\C^3$ defined by $V(xy(x+y)z)$, where $h_1=x$, $h_2=y$, $h_3=x+y$, $h_4=z$ and the $l_{ij}$ represent their pairwise intersections.   This arrangement has   $\chi_\mathcal A (t) = t^3-4t^2+5t-2$ as characteristic polynomial. 
\begin{figure}[h!]
    \centering
    \begin{tikzpicture}[scale=1.5, auto, every node/.style={circle, fill=blue!40, minimum size=4mm, inner sep=1pt}]
        \node (a) at (0.5,-1) [label=left:1] {$\mathbb C^3$};
        \node (b) at (0,0) [label=left:-1] {$h_2$};
        \node (c) at (-1,0) [label=left:-1] {$h_1$};
        \node (d) at (1,0) [label=left:-1] {$h_3$};
        \node (e) at (0,1) [label=left:1] {$l_{14}$};
        \node (f) at (2,0) [label=left:-1] {$h_4$};
	\node (g) at (1,1) [label=left:1] {$l_{24}$};
        \node (h) at (2,1) [label=left:1] {$l_{34}$};
        \node (i) at (-1,1) [label=left:2] {$l_{123}$};
        \node (j) at (0.5,2) [label=left:-2] {${\bf 0}$};

        \draw (a) -- (b);
        \draw (a) -- (c);
        \draw (a) -- (d);
        \draw (a) -- (f);
        \draw (c) -- (i);
        \draw (b) -- (i);
        \draw (d) -- (i);
        \draw (c) -- (e);
        \draw (f) -- (e);
        \draw (g) -- (b);
        \draw (g) -- (f);
        \draw (h) -- (d);
        \draw (h) -- (f);
        \draw (j) -- (i);
        \draw (j) -- (e);
        \draw (j) -- (g);
        \draw (j) -- (h);
    \end{tikzpicture}
    \caption{The arrangement $V(xy(x+y)z)$. }
    \label{fig:posetmobius2}
\end{figure}
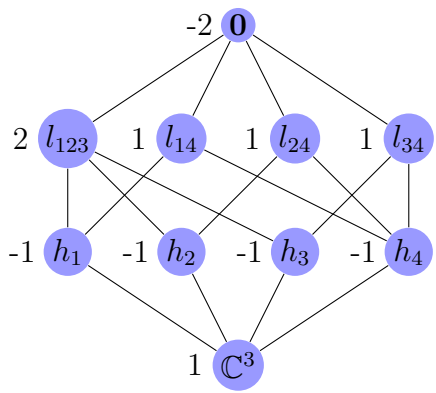

We let $w$ denote the extra variable, and let $\mathcal F$ be pencil in $\P^3$ generated by $f = xy(x+y)z$ and $w^4$.    The Gauss map $\mathcal G({\mathcal F})\colon\P^3 \dashrightarrow \P^3$ has polynomials of degree 4 as its coordinates, so we have $d_1=4$. Since we always have $d_0=1$, the only unknowns are $d_2$ and $d_3$. According to Remark~\ref{rmk:critical},   $d_3$ corresponds to  the number of critical points of the rational function $(f/w^4)|_{\P^2}$ in $U$.  A straightforward computation shows that $d_3=a_2=2$. Therefore $d_2=a_1+a_2=3+2=5.$  
\terminou
\end{example}

As  mentioned in the Introduction, it is known that the sequence $d_0, \dots, d_{n+1}$ is  {\it log-concave}, meaning that  for all $i\in\{1,\dots, n\}$, the following inequality holds
\[
d_i^2 \geq d_{i-1}d_{i+1}
\]
See \cite[Example {1.6.4}]{Laz04}.  From this and Proposition~\ref{prop:nointernal} we get the next result.  

\begin{corollary}\label{cor:chiarrang}
{Given a central affine arrangement $\mathcal{A}$ over $\C$,}
the coefficients of the characteristic polynomial $\chi_\mathcal A (t)$ form a log-concave sequence of integers with no internal zeros. 
\end{corollary}

Given a matroid $M$ which is representable over $\mathbb C$, let $\mathcal A$ be the central affine arrangement representing $M$
({see for example} \cite[Lecture 3]{Sta07}). Since  $\chi_M(t)=\chi_\mathcal A(t)$, then Corollary~\ref{cor:huhmatroid} of the introduction follows from Corollary~\ref{cor:chiarrang}. In the same spirit, we get Corollary~\ref{cor:huhgraph} via identity \eqref{eq:chiAchiG}.

\vspace{0.5cm}

\font\smallsc=cmcsc9
\font\smallsl=cmsl9

\noindent{\scriptsize\sc Universidade Federal Fluminense, Instituto de Matem\'atica e Estat\'istica.\\
Rua Alexandre Moura 8, S\~ao Domingos, 24210-200 Niter\'oi RJ,
Brazil.}

\vskip0.1cm

{\scriptsize\sl E-mail address: \small\verb?tfassarella@id.uff.br?}

{\scriptsize\sl E-mail address: \small\verb?nivaldomedeiros@id.uff.br?}


\begin{thebibliography}{123456789}

\bibitem[AHK18]{AHK18}
K. Adiprasito, J. Huh and E. Katz, 
\emph{Hodge theory for combinatorial geometries}.
 Annals of Mathematics 188 (2018), 381--452.

%
%


%
%
%
%



%
%

%

%
%
%
%


%
%
 
%
%

\bibitem[DP03]{DP03}
A. Dimca and S. Papadima,
\emph{Hypersurfaces complements, Milnor fibres and higher homotopy groups of arrangements}.
Annals of Mathematics {\bfseries 158} (2003), 473--507.




 
\bibitem[Dol24]{Dol24}
I. Dolgachev,
\emph{Classical Algebraic Geometry: a modern view I.}
Available at \url{https://dept.math.lsa.umich.edu/~idolga/CAG-1.pdf} (2024).  



%
%
%
%
%
%

\bibitem[FP07]{FP07}
T. Fassarella and J. Pereira,
\emph{On the degree of polar transformations. An approach through logarithmic foliations}.
Sel. Math. New Series {\bfseries 13} (2007), 239--252.

%
%
%

%

%
%

%
%

%
%
%

\bibitem[Her72]{Her72}
A. P. Heron, 
\emph{Matroid polynomials, Combinatorics (Proc. Conf. Combinatorial Math., Math.
Inst., Oxford, 1972)}.
Inst. of Math. and its Appl., Southend-on-Sea, (1972), 164--202.

\bibitem[Hog74]{Hog74}
S. Hoggar, 
\emph{Chromatic polynomials and logarithmic concavity}. Journal of Combinatorial Theory, 
Series B {\bfseries 16} (1974), 248--254.

\bibitem[Huh12]{Huh12}
J. Huh,
\emph{Milnor numbers of projective hypersurfaces and the chromatic polynomial of graphs}.
J. Amer. Math. Soc. (2012) 25:907--927.

%



%
%

%
%

\bibitem[Kaw78]{Kaw78}
Y. Kawamata, 
\emph{On deformations of compactifiable complex manifolds}, 
Math. Ann. 235 (1978), 247--265.

\bibitem[Laz04]{Laz04}
R. Lazarsfeld, 
\emph{Positivity in algebraic geometry}, 
vol. I and vol. II, Ergebnisseder Mathematik und ihrer Grenzgebiete. 3. Folge, 49. Springer-Verlag, Berlin, (2004).
%



%

\bibitem[Nor78]{Nor78}
Y. Norimatsu,
\emph{Kodaira Vanishing Theorem and Chern Classes for $\partial$-Manifolds}.
Proc. Japan Acad. {\bfseries 54}, Ser. A. (1978), 107--108.

%
\bibitem[OT92]{OT92}
P. Orlik and H. Terao, 
Arrangements of Hyperplanes. Grundlehren der Mathematischen Wissenschaften, 300, Springer-Verlag, Berlin, 1992.

%
%
%

\bibitem[Rea68]{Re68}
R. C. Read, 
\emph{An introduction to chromatic polynomials}. 
J. Combinatorial Theory 4 (1968), 52--71.

\bibitem[Rot70]{Rot70}
G.-C. Rota, 
\emph{Combinatorial theory, old and new, Actes du Congr\`es International des
Math\'ematiciens (Nice, 1970)}.  Tome 3, Gauthier-Villars, Paris (1971),  229--233.


%
%
%

%
%
%

\bibitem[Sil96]{Sil96}
R. Silvotti, 
\emph{On a conjecture of Varchenko}. 
Invent. Math. {\bfseries 126} (1996), 235--248.

\bibitem[Sta07]{Sta07}
R. P. Stanley,
\emph{An Introduction to Hyperplane Arrangements}.
Geometric Combinatorics, IAS/Park City Math. Ser., 13, Amer. Math. Soc., Providence, RI, 2007,  389--496.

%
%


%


\bibitem[Wel76]{Wel76}
D. Welsh, 
\emph{Matroid Theory}.
London Mathematical Society Monographs, 8, Academic Press, London-New York, (1976).

%

\end{thebibliography}
\end{document}